\pdfoutput=1
\documentclass{amsart}
\usepackage{a4wide}
\usepackage{amscd}
\usepackage{amssymb}
\usepackage{graphicx,color}
\usepackage{hyperref}
\usepackage{multirow}

\theoremstyle{plain}
\newtheorem{theorem}{Theorem}
\newtheorem{proposition}{Proposition}
\newtheorem{lemma}{Lemma}
\newtheorem{fact}{Fact}
\newtheorem{corollary}{Corollary}

\theoremstyle{definition}
\newtheorem{definition}{Definition}
\newtheorem{example}{Example}
\newtheorem{notation}{Notation}
\newtheorem{claim}{Claim}

\theoremstyle{remark}
\newtheorem{remark}{Remark}
\newcommand{\sign}{\operatorname{sign}}

\newcommand{\Dloop}
{\operatorname{\Delta^{\textup{loop}}}}

\newcommand{\rlU}{
\,~~\begin{picture}(42,8)
\put(1,3){\circle{10}}
\put(21,3){\circle{10}}
\put(41,3){\circle{10}}
\put(36.5,3){\line(-1,0){11}}
\put(16,3){\line(-1,0){11}}
\end{picture}~~
}

\newcommand{\ffone}{
\,\begin{picture}(28,8)
\put(5,3){\circle{10}}
\put(25,3){\circle{10}}
\put(20.5,3){\line(-1,0){11}}
\end{picture}~
}

\newcommand{\x}{\,
\begin{picture}(8,8)
\put(5,3){\circle{10}}
\qbezier(1,-0.5)(4.5,3)(8,6.5)
\qbezier(1.5,6.6)(5,3.3)(8.5,0)
\end{picture}~
}

\newcommand{\chordth}{
\,\begin{picture}(8,8)
\put(5,3){\circle{10}}
\qbezier(8,-1)(8,3)(8,7)
\qbezier(2,6.7)(2,3)(2,-1)
\end{picture}~
}

\newcommand{\RR}{~
\begin{picture}(48,8)
\put(5,8){\circle*{2}}
\put(25,8){\circle*{2}}
\put(45,8){\circle*{2}}
\put(5,3){\circle{10}}
\put(10,3){\vector(1,0){10}}
\put(25,3){\circle{10}}
\put(30,3){\vector(1,0){10}}
\put(45,3){\circle{10}}
\end{picture}
}

\newcommand{\LR}{
\begin{picture}(48,8)
\put(5,8){\circle*{2}}
\put(25,8){\circle*{2}}
\put(45,8){\circle*{2}}
\put(5,3){\circle{10}}
\put(20,3){\vector(-1,0){10}}
\put(25,3){\circle{10}}
\put(30,3){\vector(1,0){10}}
\put(45,3){\circle{10}}
\end{picture}
}

\newcommand{\RL}{
\begin{picture}(48,8)
\put(5,8){\circle*{2}}
\put(25,8){\circle*{2}}
\put(45,8){\circle*{2}}
\put(5,3){\circle{10}}
\put(10,3){\vector(1,0){10}}
\put(25,3){\circle{10}}
\put(40,3){\vector(-1,0){10}}
\put(45,3){\circle{10}}
\end{picture}
}
\newcommand{\LL}{~
\begin{picture}(48,8)
\put(5,8){\circle*{2}}
\put(25,8){\circle*{2}}
\put(45,8){\circle*{2}}
\put(5,3){\circle{10}}
\put(20,3){\vector(-1,0){10}}
\put(25,3){\circle{10}}
\put(40,3){\vector(-1,0){10}}
\put(45,3){\circle{10}}
\end{picture}
}

\newcommand{\LK}{
\begin{picture}(32,8)
\put(5,3){\circle{10}}
\put(10,3){\vector(1,0){11}}
\put(26,3){\circle{10}}
\end{picture}
}
\newcommand{\KL}{
\begin{picture}(32,8)
\put(5,3){\circle{10}}
\put(21,3){\vector(-1,0){11}}
\put(26,3){\circle{10}}
\end{picture}
}

\newcommand{\para}{\operatorname{parallel}}

\newcommand{\id}{\operatorname{id}}
\newcommand{\Loop}{\operatorname{loop}}

\usepackage{multirow}
\begin{document}
\title[A triple coproduct of curves and knots]{A triple coproduct of curves and knots}
\author{Noboru Ito and Takeshi Komatsuzaki}
\address{Department of Mathematics, 
Faculty of Engineering, Shinshu University, 4-17-1, Wakasato, Nagano, 380-8553, Japan
}
\email{nito@shinshu-u.ac.jp}
\address{Department of Applied Physics, School of Engineering, The University of Tokyo,  
113-8656, Japan
}
\email{st18082tk@gm.ibaraki-ct.ac.jp}
\date{December 1, 2025}
\begin{abstract}
We introduce a triple coproduct for knots on surfaces, providing a commutative framework that decomposes a single-component diagram into three components (Section~\ref{sec:preliminaries}).     This construction is motivated by the interplay between intersection theory and the affine index polynomial, and extends these ideas to a three-component setting (Section~\ref{secRelation}).  
Building on Turaev's cobracket theory, we define an integer-valued invariant under stable equivalence by combining the coproduct with an intersection-theoretic function (Theorem~\ref{curveMain}).  Unlike classical cobrackets, which often collapse distinct local configurations, our approach preserves combinatorial traces of smoothing choices, enabling fine-grained detection of local crossing patterns (Definition~\ref{LRorder}).  In the symmetric tensor setting, Reidemeister invariance uniquely determines the relations in the word space (Equations (\ref{eq:RII}), (\ref{eq:RIII})) and canonically fixes smoothing weights, revealing an intrinsic simplicity behind the algebraic framework (Corollary~\ref{cor:Unique}).  This uniqueness result positions our construction as the canonical  commutative analogue of Turaev's non-commutative cobracket and clarifies its interpretation as a classical limit of skein quantization, extending the theoretical scope beyond previously known invariants (Section~\ref{sec:Discuss}).    Examples demonstrate substantial distinguishing power, separating an infinite sequence of knots arising from distinct smoothing choices  and broadening the reach of existing invariants  (Proposition~\ref{Egprop}).   
\end{abstract} 
\maketitle
\section{Introduction}\label{intro}
Algebraic structures arising from  curves on surfaces have been well studied.   Two natural products are known: one introduced by   Goldman \cite{Goldman1986} and the other by Andersen-Mattes-Reshetikhin \cite{AndersenMattesReshetikhin1996, AndersenMattesReshetikhin1998}.   
The former yields a Lie bialgebra with the Turaev cobracket \cite{Turaev1991}; for the latter,  the Cahn operation induces the co-Jacobi and co-skew-symmetry identities \cite{Cahn2013}. 

In virtual knot theory, Kauffman \cite{Kauffman2013}, Folwaczny-Kauffman \cite{FolwacznyKauffman2013}, Cheng-Gao \cite{ChengGao2013}, and Satoh-Taniguchi \cite{SatohTaniguchi2014} independently introduced the affine index polynomial, also known as the writhe polynomial, where virtual knots are identified with stable equivalence classes of signed curves on surfaces (Turaev \cite{Turaev2006}).    Prior to these developments, Turaev \cite{Turaev2004} introduced the $u$-polynomial for virtual strings, and Henrich \cite{Henrich2010} defined a virtual knot polynomial that is related to the Goldman-Turaev Lie bialgebra.  

Motivated by these developments, we introduce a triple coproduct that unifies ideas from intersection theory and skein-theoretic constructions, providing a new perspective on invariants for knots on surfaces.  Let $\mathcal{D}$ be the set of stable homeomorphism classes of oriented knot diagrams on oriented surfaces; let $\mathcal{C}^{3}$ be the set of stable homeomorphism classes of three-component oriented curves  on oriented surfaces.  Either $D$ or a tuple $\left(C^{(1)}, C^{(2)}, C^{(3)}\right)$ denotes a stable homeomorphism class of knot diagrams or of three-component curves, respectively.             
\begin{figure}[h]
\includegraphics[width=5cm]{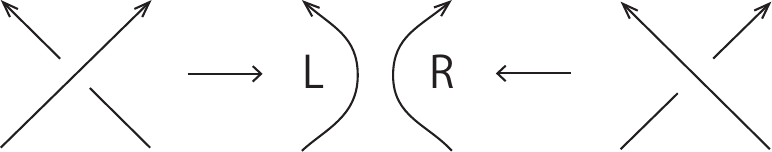}
\caption{Smoothing.  The label $L$ (resp.~$R$) indicates ``left'' (resp.~``right'').}\label{Seifert}
\end{figure}  
For $D \in \mathcal{D}$, an unordered pair $a, b$ of crossings in $D$ is called \emph{parallel}  if (and only if) smoothing these two crossings as in Figure~\ref{Seifert} produces a three-component curve $C^{(1)}_{ab}, C^{(2)}_{ab}, C^{(3)}_{ab}$.   For each parallel pair $a, b$, we assign an element $w_{ab}$ of the word space $W_{LRNP}$ (Definition~\ref{EqRL}) and we define a map $\Delta$ by    
\begin{equation}
\label{formula:coproduct}
\Delta : \mathcal{D} \to \mathbb{Z}[\mathcal{D}^3]  \otimes W_{LRNP};\quad D \mapsto \sum_{a, b : \para}  \varepsilon_a  \varepsilon_b \left(C^{(1)}_{ab}, C^{(2)}_{ab}, C^{(3)}_{ab}\right) \otimes w_{ab}, 
\end{equation}    
where $\varepsilon_a$ (resp.~$\varepsilon_b$) denotes the local writhe of $a$ (resp.~$b$).   
This $\Delta$ does not preserve  stable equivalence, but $(\nu \otimes \id) \circ \Delta$ does; $\nu$ is an extension of the intersection number to three-component curves   (Definition~\ref{nu}).   The value $(\nu \otimes \id) \circ \Delta(D)$ lies in $W_{LRNP}$.  
\begin{theorem}\label{curveMain}
Let $D$ be a stable homeomorphism class of a knot diagram.  Then $(\nu \otimes \id) \circ \Delta(D)$ is invariant under stable equivalence.   
\end{theorem}
Further, Section~\ref{sec:Discuss} shows that, when knot diagrams are mapped into the symmetric tensor setting following Turaev's skein-theoretic framework, Reidemeister invariance imposes necessary and sufficient relations on the word space.  This condition forces a unique normalization of smoothing weights, yielding a canonical extension  
\[
\Dloop : \mathcal{D} \to \textup{Sym}^3 (V)  \otimes W_{LRNP}, 
\]
and proving that our construction is the unique extension compatible with Turaev's skein-theoretic framework.    
This result clarifies the theoretical  position of our construction as the canonical commutative analogue of Turaev's non-commutative cobracket and its interpretation as a classical limit of skein quantization.     
\section{Preliminaries}
\label{sec:preliminaries} 
\begin{figure}[h]
\includegraphics[width=5cm]{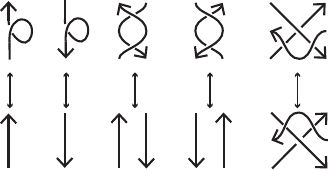}
\caption{A generating set of oriented Reidemeister moves: ${\Omega_{1a}}$, ${\Omega_{1b}}$, ${\Omega_{2c}}$, ${\Omega_{2d}}$, and ${\Omega_{3a}}$ from the left to the right.}
\label{Deformations}
\end{figure}

\begin{figure}[htbp] 
   \centering
   \includegraphics[width=2in]{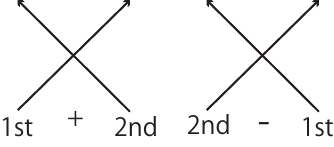} 
   \caption{A self-intersection of link projection is shown on the left.  If we pass through the self-intersection from  lower left to upper right first, the intersection number is $+1$ (center); otherwise, it is $-1$ (right).} 
   \label{fig:insec}
\end{figure}

\begin{figure}[h!]
   \includegraphics[width=5cm]{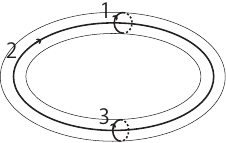} 
 \caption{A three component curve}
   \label{nuEgA}
\end{figure}

\begin{definition}[link diagram, stable homeomorphism]
A \emph{link diagram} is the image of a generic immersion of oriented circles into an oriented closed surface, where each self-intersection is assigned over/under information and is called a \emph{crossing}.     
If the over/under information of the crossings of a link diagram is ignored, the result  is called a curve or a link projection.  Each self-intersection of a curve is called a crossing if there is no confusion or called a double point when we explicitly wish to regard it as \emph{not} a crossing.     
Two link diagrams are said to be \emph{stably homeomorphic} if there exists  a homeomorphism of their regular neighborhoods in the ambient surfaces that maps  the first link diagram onto the second one and preserves the orientations of both the diagram and the surface.    
\end{definition}
\begin{remark}
Adding handles to the ambient surface away from the neighborhood of the link diagram, does not change its  stable homeomorphism class.    
\end{remark}
\begin{definition}[stable equivalence]
Two one-component link diagrams are said to be \emph{stably equivalent} if they can be related by finite sequence of deformations $\Omega_{1a}$, ${\Omega_{1b}}$, ${\Omega_{2c}}$, ${\Omega_{2d}}$, and ${\Omega_{3a}}$ shown in Figure~\ref{Deformations}, up to stable homeomorphisms \footnote{The set of moves $\Omega_{1a}$, ${\Omega_{1b}}$, ${\Omega_{2c}}$, ${\Omega_{2d}}$, and ${\Omega_{3a}}$ forms a generating set for oriented Reidemeister moves (see \cite{ItoIwamoto2025}).}.    
\end{definition}
\begin{definition}[chord diagram]
A configuration of $n$ pairs of points on circles, considered up to ambient isotopy and reflection of the circles, is called a \emph{chord diagram}.  For each pair, two points of are traditionally connected by a straight arc, called a \emph{chord}.     
In particular, if an immersion $f : \to \mathbb{R}^2$ represents a link diagram $D$, each crossing $c$ of $D$ has a preimage $f^{-1}(c)$ consisting of exactly two points, which form a pair on $S^1$.  Thus $f^{-1}(D)$ specifies a chord diagram, which we call the chord diagram of $D$.   
\end{definition}  
\begin{fact}\label{fact:ch}
Every link diagram obtains a chord diagram.  
\end{fact}
\begin{remark}\label{RMKpara}
If the reader is familiar with chord diagrams, a \emph{parallel} pair is symbolically represented by $\chordth$; while a non-parallel pair is represented by $\x$.     
\end{remark}
\begin{definition}
[$\Delta$ and words of $L, R, N, P$]\label{LRorder}
For any crossing of a link diagram $D$ on a surface, we regard the path entering from the lower right as positive (``$+$'') and the path entering from the lower left as negative (``$-$'').   In the chord diagram corresponding to $D$, each of the two ends of every chord is then labelled with the associated sign $\pm$.  
Consequently, each endpoint of each chord is assigned a unique sign.  
We orient the circle in the chord diagram counterclockwise (we fix the circle orientation of any chord diagram in this way throughout this paper).  Then we assign the label $R$ just before a ``$+$'' endpoint is encountered, and $L$ just after it is encountered along the circle orientation.  Note that this convention corresponds to the left ($L$) and right ($R$) sides in the smoothing shown in Figure~\ref{Seifert}.    
Therefore, for every chord, the meanings of the labels $L$ and $R$ are unambiguously determined. 
In general, each crossing $a$ of a link diagram has a local writhe sign $\varepsilon_a$.  To distinguish the local writhe sign from the signs assigned to chord endpoints above, we consistently denote it by $P$ (positive) or $N$ (negative) throughout this paper.  

Using the above data for two crossings $a$ and $b$ to be smoothed, the patterns of parallel pairs---each  represented by two parallel chords $\chordth$---fall into  ten possible cases (Figure~\ref{OneToTheree}), from which we derive the LRNP word $w_{a, b}$.   Figure~\ref{sProcess} illustrates how a parallel pair yields a word by applying the two smoothings.  
Note that a rotating the chord diagram by $\pi$ in the plane yields the relations $LLNP=LLPN$, $RRNP=RRPN$, and $LR\ast\ast=RL\ast\ast$ ($\ast=N$ or $P$).  

Next, we smooth two crossings $a, b$ that correspond to a parallel pair.  Then we choose arbitrary ordering of the resulting components, denoting the first, second, and third components by $C^{(1)}_{ab}$, $C^{(2)}_{ab}$, $C^{(3)}_{ab}$ (but it will later show that our invariant is independent of the chosen order of the three components).  

Let  $\mathbb{Z}[\mathcal{D}^3]$ denote the $\mathbb{Z}$-module  generated by $\mathcal{D}^3$.  Then we define 
\begin{equation}
\Delta(D) := \sum_{a, b : \para} \varepsilon_a \varepsilon_b (C^{(1)}_{ab}, C^{(2)}_{ab}, C^{(3)}_{ab}) \otimes w_{ab},   
\end{equation}
where the sum runs over the unordered parallel pairs of crossings.  This leads to the definition of our invariant as follows: 
\begin{equation}
(\nu \otimes \id) \Delta(D) := \sum_{a, b : \para} \varepsilon_a  \varepsilon_b \nu(C^{(1)}_{ab}, C^{(2)}_{ab}, C^{(3)}_{ab}) \otimes w_{ab}.  
\end{equation}
\
\begin{figure}
\includegraphics[width=14cm]{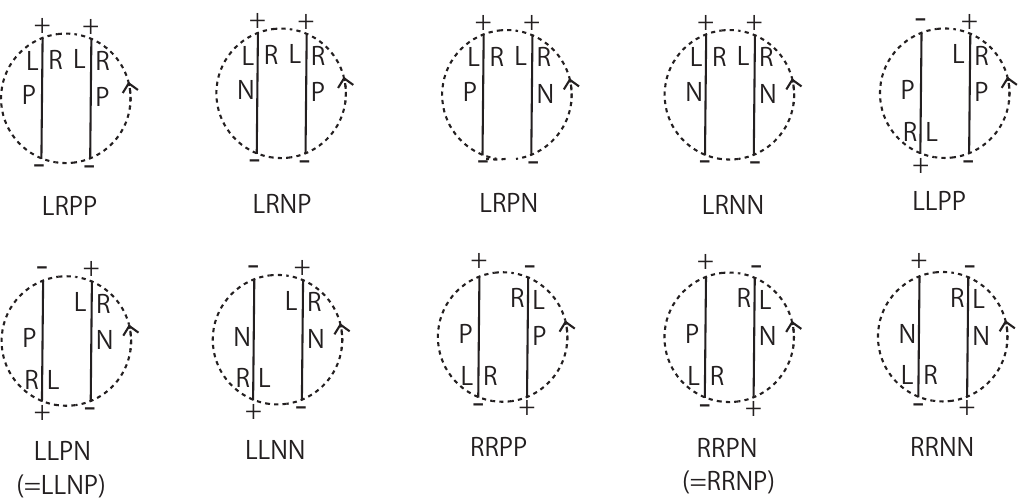}
\caption{Ten possible cases corresponding to words of L, R, N,  P.}\label{OneToTheree}
\end{figure}              
\end{definition}
\begin{figure}[htbp]  
   \includegraphics[width=12cm]{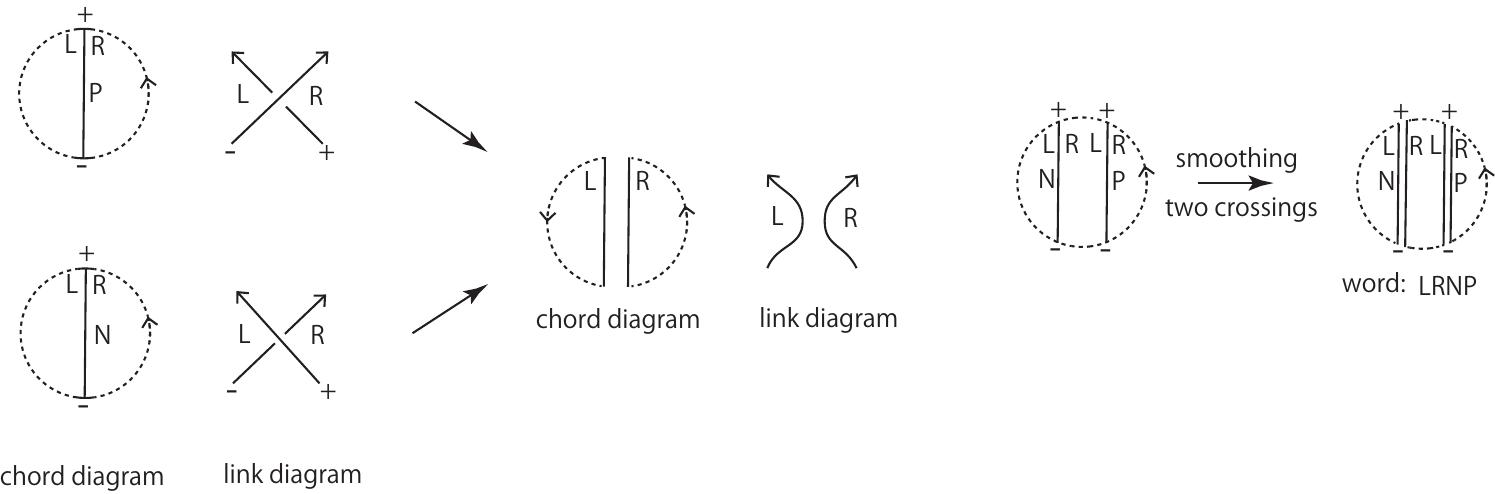}
\caption{
An example illustrating the  procedure for handling a parallel pair of crossings.  By definition, a parallel pair consists of two crossings whose chords in the chord diagram are parallel.  We start from the chord diagram and its corresponding link diagram, and  smooth the two crossings according to the convention shown in Figure~\ref{Seifert}.    The resulting three-component chord diagram contains a distinguished center component.  By going along this center component one counterclockwise around the circle, we read off the labels on the chord endpoints in order.  Specifically, we first write the two-letter LR word corresponding to the parallel pair, followed by the two-letter NP word corresponding to the parallel pair.  This determines the word in $\{L, R, N, P \}$, which encodes one of the ten possible cases listed in Figure~\ref{OneToTheree} (e.g., LRNP, as shown in this example).      
}\label{sProcess}
\end{figure}
\begin{definition}[invariant $\nu$]\label{nu}
Every link diagram is regarded as a link projection.  For each link projection $D$, we use the same symbol $D$ to denote the link projection obtained by ignoring over/under information of its  crossings.  
Let $I(D, D')$ denote the intersection number of $D$ and $D'$ in the same surface.       
Let $\left(C^{(i)}, C^{(j)}, C^{(k)}\right)$ be a three-component link projection with an order $(i, j, k)$.  
For a $3$-component link projection $C=\left(C^{(1)}, C^{(2)}, C^{(3)}\right)$, let 
\begin{align*}
\nu (C) &:= I(C^{(1)}, C^{(2)})  I(C^{(1)}, C^{(3)}) + I(C^{(2)}, C^{(1)})  I(C^{(2)}, C^{(3)}) + I(C^{(3)}, C^{(1)})  I(C^{(3)}, C^{(2)}).    
\end{align*}     
\end{definition}
By definition we have
\begin{proposition}\label{PropNu}Let $C$ be a three-component curve (i.e., link projection) on a surface.  The integer-valued function $\nu (C)$ is invariant under stable equivalence of link projections.  
\end{proposition}
\begin{notation}
We record a convention here. The intersection number $I(C^{(i)}, C^{(j)})$ of two link diagrams on surfaces $C^{(i)}$ and $C^{(j)}$ is computed from their underlying link  \emph{projections}, which are generic curves on surfaces,  using the standard definition of intersection number (Figure~\ref{fig:insec}).   
\end{notation}
\begin{example}
For a curve in Figure~\ref{nuEgA}, 
\begin{align*}
\nu \left(C^{(1)}, C^{(2)}, C^{(3)}\right) &= I(C^{(1)}, C^{(2)})  I(C^{(1)}, C^{(3)}) + I(C^{(2)}, C^{(1)})  I(C^{(2)}, C^{(3)}) + I(C^{(3)}, C^{(1)})  I(C^{(3)}, C^{(2)}) \\
&= 0+I(C^{(2)}, C^{(1)})  I(C^{(2)}, C^{(3)})+0 = (+1) \cdot (-1) = -1.  
\end{align*}
Here, the intersection numbers  $I(C^{(2)}, C^{(1)})$ and $I(C^{(2)}, C^{(3)})$ are given as follows. 
\begin{itemize}
\item $I(C^{(2)}, C^{(1)})$. 
The number $I(C^{(2)}, C^{(1)})$ is computed according to the order of arguments: the first (resp.~second) argument is treated as the first (resp.~second) component, regardless of each index.  Hence, for $I(C^{(2)}, C^{(1)})$, $C^{(2)}$ is the \emph{first} component and $C^{(1)}$ is the \emph{second} component, which corresponds to the pattern as in the left of Figure~\ref{fig:insec}, giving $+1$.    
\item The number $I(C^{(2)}, C^{(3)})$ is computed with the first component $C^{(2)}$ by the first argument and the second component $C^{(3)}$ by the second argument, which corresponds to the pattern as in the right of Figure~\ref{fig:insec}, giving $-1$.
 \end{itemize}
\end{example}
\begin{definition}[$W_{LRNP}$]\label{EqRL}
Let $W_{LRNP}$ denotes the quotient $\mathbb{Z}$-module generated by the set of words of length four with  following generators and relations: 

\noindent Generators: $LLPP$, $LLNP$ (it is permitted to represent it by $LLPN$), $LLNN$, $LRPP$, $LRNP$, $LRPN$, $LRNN$, $RRPP$, $RRNP$ (it is permitted to represent it by $RRPN$), $RRNN$.  

\noindent Relations: 
\begin{align}
LLPP=LRNP=RRNN,& \ LLNN=LRPN=RRPP,  \  LLPN=LRPP=LRNN=RRPN, \label{eq:RII}\\
RRPP=2RRPN, & LLPP=2LLPN. \label{eq:RIII}
\end{align}   
\end{definition}
\begin{remark}
The generators correspond to ten possible parallel pairs.  The relation (\ref{eq:RII}) corresponds to the invariance of the second Reidemeister moves $\Omega_{2c}$ and $\Omega_{2d}$. The relation (\ref{eq:RIII}) corresponds to the invariance of the third Reidemeister move $\Omega_{3a}$.   
\end{remark}
\begin{notation}[$\Delta|_{ab}$]\label{not:state}  
Let $\Delta|_{ab} (D)$ $=$ $\varepsilon_a \varepsilon_b  \left(C^{(1)}_{ab}, C^{(2)}_{ab}, C^{(3)}_{ab}\right) \otimes w_{ab}$, then $\displaystyle \Delta(D) = \sum_{a, b : \para} \Delta|_{ab} (D)$.  
We allow the use of the symbol $\nu \Delta$ (resp.~$\nu \Delta|_{ab}$) to indicate the composition $(\nu \otimes \id) \circ \Delta$ (resp.~$(\nu \otimes \id) \circ \Delta|_{ab}$) when there is no risk of confusion.
\end{notation}
\section{Proof of Theorem~\ref{curveMain}}\label{sec:ProofM}
We are now prepared to prove Theorem~\ref{curveMain}.        
\subsection{Invariance of ${\Omega_{1a}}$ and ${\Omega_{1b}}$}\label{secRI} 
Suppose that the move ${\Omega_{1a}}$ (or ${\Omega_{1b}}$) on a diagram $D$ creates a single crossing $A$, and let $D'$ be the resulting diagram.  Then we have 
\begin{equation}\label{eq:Rone}
\nu \Delta (D') - \nu \Delta (D) = \sum_{\substack{a, b : \para \\ a=A~{\textrm{or}~} b=A}} \varepsilon_a \varepsilon_b \nu  \left(C^{(1)}_{ab}, C^{(2)}_{ab}, C^{(3)}_{ab}\right) \otimes w_{ab}.   
\end{equation}
For a parallel pair $(A, b)$ or $(a, A)$ in $D$, since there is a component $C^{(k)}_{ab}$ for some $k \in \{1, 2, 3\}$ with no crossings, we have 
\begin{align}\label{RIzero}
\nu \left(C^{(1)}_{ab}, C^{(2)}_{ab}, C^{(3)}_{ab} \right) = 0,   
\end{align}
which implies the invariance of $\nu \Delta$ under the first Reidemeister moves; more precisely, the right-hand side of (\ref{eq:Rone}) is zero, i.e.,    
\[\nu \Delta(D')=
\sum_{\substack{a, b : \para \\ a \neq A~{\textrm{and}~} b \neq A}} \varepsilon_a \varepsilon_b \nu  \left(C^{(1)}_{ab}, C^{(2)}_{ab}, C^{(3)}_{ab}\right) \otimes w_{ab} = \nu \Delta(D).  
\]
\subsection{Invariance of ${\Omega_{2c}}$ and ${\Omega_{2d}}$}\label{secRII}
\begin{figure}
\includegraphics[width=13cm]{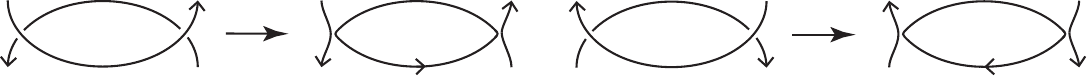}
\caption{Bigon with crossings $A, B$ applied by ${\Omega_{2c}}$ (left), ${\Omega_{2d}}$ (right).      
}\label{RII}
\end{figure}
\begin{figure}[htbp] 
   \centering
   \includegraphics[width=8cm]{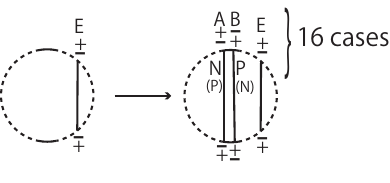}
   \caption{The 16 cases corresponding to the configuration  represented by chord diagrams}
   \label{InvRII}
\end{figure}

Suppose that the move ${\Omega_{2c}}$ or ${\Omega_{2d}}$ on a link diagram $D$ generates two crossings $A$ and $B$, and let $D'$ be the resulting diagram.  
Firstly, Figure~\ref{RII} implies the following.   
\begin{lemma}\label{RIIpara}
\[
\nu \left( C^{(1)}_{AB}, C^{(2)}_{AB}, C^{(3)}_{AB} \right) = 0.  
\] 
\end{lemma}      
Hence we have 
\begin{align}\label{eqRIIdiff}
\nu \Delta (D') - \nu \Delta (D) &= \sum_{\substack{a, b : \para \\ 
|\{ a, b \} \cap \{ A, B \}  |=1}} \varepsilon_a \varepsilon_b \nu \left( C^{(1)}_{ab}, C^{(2)}_{ab}, C^{(3)}_{ab} \right) \otimes w_{ab}.  
\end{align}
We now check the cases where $ 
|\{ a, b \} \cap \{ A, B \}  |=1$ for any parallel pair $a, b$.    
Let $E$ be a crossing.  We will consider the situations such that pairs $(E, A)$, $(A, E)$, $(E, B)$, or $(B, E)$ are parallel; their configuration is represented by chord diagrams as in Figure~\ref{InvRII}.
The 16 cases are listed as in Table~\ref{RIICases}.  Since $\varepsilon_{A} = - \varepsilon_{B}$, the contribution to the whole sum $\pm(w_{EB}-w_{EA})$ in each case.  However, $w_{EB}-w_{EA}=0$ holds by relations of the definition $W_{LRNP}$ as in  Definition~\ref{EqRL}, which implies the invariance under $\Omega_{2c}$ and $\Omega_{2d}$.  
\begin{table}[h]
\label{RIICases}
\begin{tabular}{|c|c|c|c|c|}\hline
crossing $A$ & crossing $B$ & crossing $E$ & $w_{E, B}$ by smoothing $E, B$ & $w_{E, A}$ by smoothing $E, A$\\ \hline
$-$, $N$&$+$, $P$&$+$, $P$&$LRPP$& $LLNP$ \\ \hline
$-$, $N$&$+$, $P$&$+$, $N$&$LRPN$& $LLNN$ \\ \hline
$-$, $P$&$+$, $N$&$+$, $P$&$LRNP$& $LLPP$  \\ \hline
$-$, $P$&$+$, $N$&$+$, $N$&$LRNN$& $LLPN$ \\ \hline
$+$, $N$&$-$, $P$&$+$, $P$&$LLPP$& $LRNP$ \\ \hline
$+$, $N$&$-$, $P$&$+$, $N$&$LLPN$& $LRNN$ \\ \hline
$+$, $P$&$-$, $N$&$+$, $P$&$LLNP$& $LRPP$ \\ \hline
$+$, $P$&$-$, $N$&$+$, $N$&$LLNN$& $LRPN$ \\ \hline
$-$, $N$&$+$, $N$&$-$, $P$&$RRPP$& $LRPN$ \\ \hline
$-$, $N$&$+$, $N$&$-$, $N$&$RRPN$& $LRNN$ \\ \hline
$-$, $P$&$+$, $N$&$-$, $P$&$RRNP$& $LRPP$ \\ \hline
$-$, $P$&$+$, $N$&$-$, $N$&$RRNN$& $LRNP$ \\ \hline
$+$, $N$&$-$, $P$&$-$, $P$&$LRPP$& $RRNP$ \\ \hline
$+$, $N$&$-$, $P$&$-$, $N$&$LRNP$& $RRNN$ \\ \hline
$+$, $P$&$-$, $N$&$-$, $P$&$LRPN$& $RRPP$ \\ \hline
$+$, $P$&$-$, $N$&$-$, $N$&$LRNN$& $RRPN$ \\ \hline
\end{tabular}
\end{table}
$\hfill \Box$ 
\subsection{Invariance of ${\Omega_{3a}}$}\label{secRIII}
Let $A$, $B$, and $C$ be crossings of the triangle involved in the move ${\Omega_{3a}}$ as in Figure~\ref{Deformations}.   
\begin{figure}[htbp] 
   \centering
   \includegraphics[width=10cm]{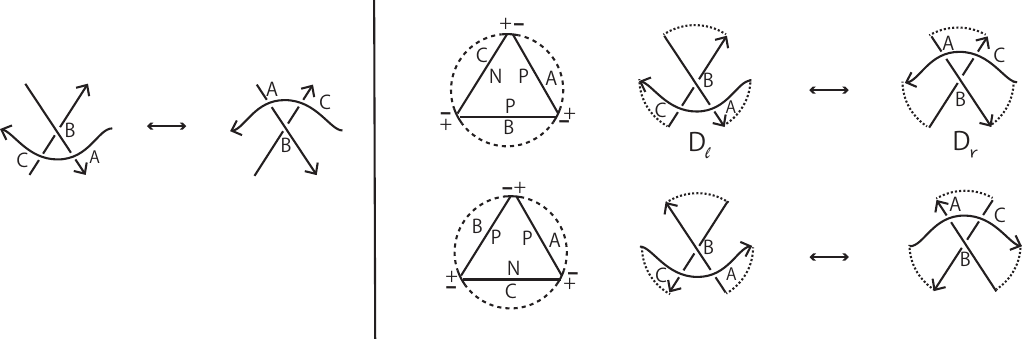} 
   \caption{A single ${\Omega_{3a}}$ (left half) is among two oriented diagrams $D_{\ell}$ and $D_r$, each of which has two cases (right half) with chord diagrams of $D_{\ell}$.}
   \label{IIIaCases}
\end{figure}      
Let $D_r$ and $D_{\ell}$ be the diagrams as in  Figure~\ref{IIIaCases}.   For each $\bullet \in \{\ell, r\}$, define  
\begin{align}\label{eqRIIIdiff} 
\sum \ast_k (D_{\bullet}) &= \sum_{\substack{a, b : \para~{\textrm{in}}~D_{\bullet} \\ 
 |\{ a, b \} \cap \{ A, B, C \}  |= k}} \varepsilon_a \varepsilon_b \nu \left( C^{(1)}_{ab}, C^{(2)}_{ab}, C^{(3)}_{ab}  \right) \otimes w_{ab}.
 \end{align}  
\subsubsection{\bf Pair $a, b$ containing exactly \underline{zero} or \underline{one} element of $\{ A, B, C \}$}\label{secZero}
\begin{figure}[h!]
\includegraphics[width=12cm]{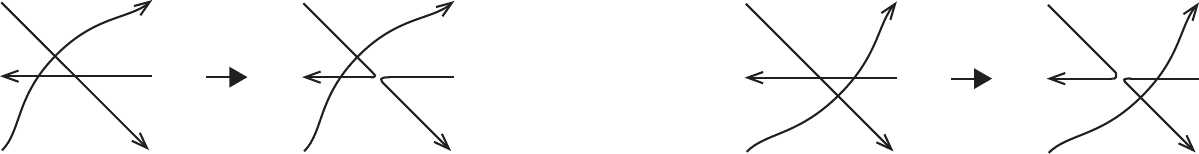} 
\caption{A pair of triangles  which will be applied by ${\Omega_{3a}}$ and one crossing of each triangle is smoothed.        
}\label{RIIIaDiag}
\end{figure}

By definition, we have $\sum \ast_0 (D_{\ell})$ $=$ $\sum \ast_0 (D_r)$.  By Figure~\ref{RIIIaDiag} and using the invariance of $\nu$ under stable equivalence (Proposition~\ref{PropNu}), we also have $\sum \ast_1 (D_{\ell})$ $=$ $\sum \ast_1 (D_r)$.   
\subsubsection{\bf Pair $a, b$ containing exactly \underline{two} elements of $\{ A, B, C \}$}\label{secTwo}
Seeing $D_{\ell}$ (resp.~$D_r$) in   Figure~\ref{IIIaCases}, any pair of crossings in  $\{ A, B, C \}$ is parallel (resp.~\emph{not} parallel).  Hence it suffices to prove the following (cf.~Notation~\ref{not:state}).   
\begin{lemma}\label{SumZero} 
\begin{align}\label{eqSumZero}
\nu \Delta|_{AB} (D_{\ell})  + \nu \Delta|_{BC} (D_{\ell})  + \nu \Delta|_{AC} (D_{\ell}) = 0.
\end{align}  
\end{lemma}
\noindent{\bf Proof of Lemma~\ref{SumZero}.}  
First, considering the upper case in Figure~\ref{IIIaCases},  
\begin{align*}
\nu \Delta|_{AB} (D_{\ell}) &= \nu  {\left(C^{(1)}_{AB}, C^{(2)}_{AB}, C^{(3)}_{AB} \right)} \otimes RRPP, 
\nu \Delta|_{AC} (D_{\ell}) &= - \nu {\left(C^{(1)}_{AC}, C^{(2)}_{AC}, C^{(3)}_{AC} \right)} \otimes RRPN,\\~{\textrm{and}}~
\nu \Delta|_{BC} (D_{\ell}) &= -  \nu {\left(C^{(1)}_{BC}, C^{(2)}_{BC}, C^{(3)}_{BC} \right) \otimes RRNP}.  
\end{align*}
Noting that $\nu  {\left(C^{(1)}_{AB}, C^{(2)}_{AB}, C^{(3)}_{AB} \right)} = \nu  {\left(C^{(1)}_{AC}, C^{(2)}_{AC}, C^{(3)}_{AC} \right)} = \nu  {\left(C^{(1)}_{BC}, C^{(2)}_{BC}, C^{(3)}_{BC} \right)}$, 
\[
\nu \Delta|_{AB} (D_{\ell}) + \nu \Delta|_{AC} (D_{\ell}) + \nu \Delta|_{BC} (D_{\ell}) = \nu  {\left(C^{(1)}_{AB}, C^{(2)}_{AB}, C^{(3)}_{AB}  \right)} \otimes (RRPP - RRPN - RRNP). 
\]

Second, considering the lower case in Figure~\ref{IIIaCases}, 
\begin{align*}
\nu \Delta|_{AB} (D_{\ell}) &= \nu  {\left(C^{(1)}_{AB}, C^{(2)}_{AB}, C^{(3)}_{AB} \right)} \otimes LLPP, 
\nu \Delta|_{AC} (D_{\ell}) &= - \nu {\left(C^{(1)}_{AC}, C^{(2)}_{AC}, C^{(3)}_{AC} \right)} \otimes LLPN,\\~{\textrm{and}}~
\nu \Delta|_{BC} (D_{\ell}) &= -  \nu {\left(C^{(1)}_{BC}, C^{(2)}_{BC}, C^{(3)}_{BC} \right) \otimes LLNP}.  
\end{align*}
Noting that $\nu  {\left(C^{(1)}_{AB}, C^{(2)}_{AB}, C^{(3)}_{AB} \right)} = \nu  {\left(C^{(1)}_{AC}, C^{(2)}_{AC}, C^{(3)}_{AC} \right)} = \nu  {\left(C^{(1)}_{BC}, C^{(2)}_{BC}, C^{(3)}_{BC} \right)}$, 
\[
\nu \Delta|_{AB} (D_{\ell}) + \nu \Delta|_{AC} (D_{\ell}) + \nu \Delta|_{BC} (D_{\ell}) = \nu  {\left(C^{(1)}_{AB}, C^{(2)}_{AB}, C^{(3)}_{AB}  \right)} \otimes (LLPP - LLPN - LLNP). 
\]
This together with Definition~\ref{EqRL} implies the formula (\ref{eqSumZero}) holds.  
\hfill({\bf End of Proof of Lemma~\ref{SumZero}}.) $\Box$

\noindent{\bf Proof of Theorem~\ref{curveMain}.}
Combining the above results,  
\begin{align*}
\nu \Delta (D_{\ell}) &= \sum \ast_0 (D_{\ell}) +  \sum \ast_1 (D_{\ell}) + \sum \ast_2 (D_{\ell})  \\
&\stackrel{\textrm{Lemma~\ref{SumZero}}}{=} \sum \ast_0 (D_{\ell}) +  \sum \ast_1 (D_{\ell}) \\
&\stackrel{\textrm{Section~\ref{secZero}}}{=}  \sum \ast_0 (D_r) +  \sum \ast_1 (D_r) 
= \nu \Delta (D_r).  
\end{align*}
\hfill({\bf End of Proof of Theorem~\ref{curveMain}}.) $\Box$
\section{Examples}   
\begin{proposition}
\label{Egprop}
There exist infinitely many distinct  links $\{ L_n \}$ on surfaces as in  Figure~\ref{EgCurve} such that if $i \neq j$, $L_i$ and $L_j$ are not stably equivalent.  
\begin{figure}[h] 
   \centering
\includegraphics[width=14cm]{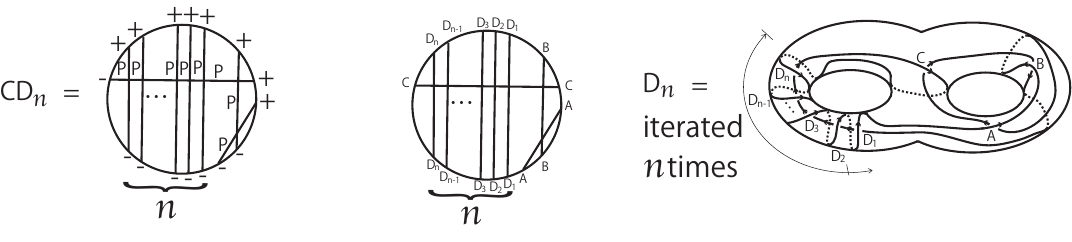} 
   \caption{$L_n$ corresponding to the number of $n$ chords ($n=1, 2, 3, \dots$) appearing in the position in the chord diagram (left) and $L_n$ on a surface (right).}
   \label{EgCurve}
\end{figure}
\end{proposition}
\begin{figure}[htbp] 
   \centering
   \includegraphics[width=14cm]{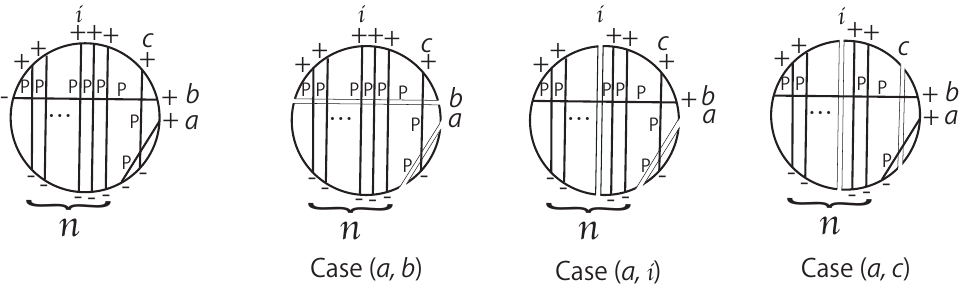} 
   \caption{Each contribution to each term.}
   \label{nCurvePcal}
\end{figure}
\begin{proof}
Let $n$ be a given positive integer.  We consider the link diagram on a surface (Figure~\ref{EgCurve}) corresponding to the chord diagram $CD_n$ as in Figure~\ref{nCurvePcal}.    Let $(a, b)$, $(a, i)$, $(c, i)$ ($i=1,2, \dots, n$) be the parallel pairs in $CD_n$ where, written as ordered pairs $(a, b)$ for convenience, although the order is irrelevant.  Then we have 
\begin{align*}
(\nu \otimes \id) \Delta(D) &= \sum_{a, b : \para} \varepsilon_a  \varepsilon_b \nu(C^{(1)}_{ab}, C^{(2)}_{ab}, C^{(3)}_{ab}) \otimes w_{ab}\\
&= n LRPP + \sum_{i=1}^n LRPP + \sum_{i=1}^n LRPP,  \\
\end{align*}
where the first term corresponds to the parallel pair $(a, b)$, the second term $(a, i)$, and the third term to $(c, i)$  ($i=1,2,\cdots, n$).  

Hence, for curves with different $n$, the values $(\nu \otimes \id) \Delta(D)=3nLRPP$ are distinct.   
\end{proof}
\section{The relationship with the affine index polynomial}\label{secRelation}
Throughout this section, to avoid confusion, if self-intersections appear in (virtual) knot diagrams, we call them \emph{crossings}; if self-intersections appear in diagrams of curves on surfaces, we call them \emph{double points}. 

A knot diagram on a surface corresponds to a curve with over/under information at double points; further, the stable equivalence classes of curves on surfaces are called \emph{virtual strings}, and an element is often represented by a flat virtual knot diagram with a base point.  
    
We recall the construction of the affine index polynomial \cite{HigaNakamuraNakanishiSatoh2022}.    
Let $\gamma_i$ and $\bar{\gamma}_i$ be the components of a two-component virtual link diagram obtained by smoothing the $i$th crossing, where $\gamma_i$ (resp.~$\bar{\gamma}_i$) is the right (resp.~left) component $R$ (resp.~$L$) appearing in Figure~\ref{Seifert}.  Let $I(\gamma_i, \bar{\gamma}_i)$ be the intersection number of two curves $\gamma_i$ and  $\bar{\gamma}_i$.  \footnote{$\gamma_i \cdot  \bar{\gamma}_i$ indicates the intersection number in \cite{HigaNakamuraNakanishiSatoh2022}.}.   
\begin{theorem}[The affine index polynomial \cite{ChengGao2013, FolwacznyKauffman2013, Kauffman2013, SatohTaniguchi2014}]\label{affinethm}
For a virtual knot $K$, let $D$ be an $n$-crossing knot diagram, $c_1, c_2, \dots, c_n$ crossings, and $\varepsilon_i$ the local writhe of $c_i$.  
The Laurent polynomial 
\[ W_D (t) = \sum_{i=1}^{n} \varepsilon_i (t^{I(\gamma_i, \bar{\gamma}_i)} - 1) \]
is an invariant of $K$.     
\end{theorem}
Examining the polynomial $W_D(t)$, we observe that 
\[
\frac{d}{dt} W_D (t)|_{t=1} = \sum_{i=1}^{n} \varepsilon_i I(\gamma_i, \bar{\gamma}_i).  
\]
In order to facilitate comparison with Theorem~\ref{affinethm}, we use the symbol $\left(C^{(1)}_i, C^{(2)}_i \right)$ to denote  the pair of curves $(\gamma_i, \bar{\gamma}_i)$.   
The intersection number can also be  expressed as a bilinear function, known as a Gauss diagram formula, $\langle \ffone, G_C \rangle$, for a Gauss diagram $G_C$ of a curve $C$ on a surface.   Hence 
\[
\left \langle  \ffone, \left(C^{(1)}_i, C^{(2)}_i \right) \right \rangle = I(\gamma_i, \bar{\gamma}_i).  
\] 
We summarize the comparison below.  
\begin{center}
\begin{tabular}{|c|c|c|}\hline
Invariants & $\frac{d}{dt} W_D (t)|_{t=1}$ & $(\nu \otimes \id) \circ \Delta$   \\ \hline
Transit objects & $\left(C^{(1)}_i, C^{(2)}_i \right)$ & $\left(C^{(1)}_{ab}, C^{(2)}_{ab}, C^{(3)}_{ab}\right) \otimes w_{ab}$ \\ \hline 
Smoothing & a single crossing & two crossings \\ \hline
Gauss diagrams  & \ffone & \rlU \\\hline
Intersection number  & $I \left(C^{(1)}_i, C^{(2)}_i \right)$ & 
$\nu$ \\\hline
Sum  & $\displaystyle \sum_i \sign (c_i)\ \cdot$ (term) & $\displaystyle \sum_{a, b: \para} \varepsilon_a, \varepsilon_b$ (term) $\otimes$ (word)  \\\hline
\end{tabular}
\end{center}
We also note that the linking number $lk(L)$ is $\langle \frac{1}{2}(\LK+\KL), G_L \rangle$ \cite{PolyakViro1994}.  Analogously, Milnor's triple linking number $\mu(L)$ is a linear combination of $\LL$~, $\LR$~, $\RL$~, and $\RR$~.       

\begin{remark}
Turaev introduced an operation \footnote{Turaev use the symbol $\nu$, but we shall refer to it here as $\nu_T$ to avoid confusion.}, which yields cobracket.  
One may be curious about comparing the second iteration $\nu_T^{(2)}$ $=$ $(\id \otimes \nu_T) \nu_T$ of Turaev's virtual  cobracket $\nu_T$ with our operation.    Turaev's operation, for a homotopy class $\langle \alpha \rangle$ of  virtual string $\alpha$, is defined as 
\begin{align}\label{Tcb}
\nu_T (\langle \alpha \rangle )&\sum_{e : {\text{arrow}}} = \langle \alpha^1_e \rangle \otimes  \langle \alpha^2_e \rangle - \langle \alpha^2_e \rangle\otimes \langle \alpha^1_e \rangle, 
\end{align}
where an arrow $e$ corresponds to a crossing which will be smoothed (for the precise definition, see \cite{Turaev2004}).   
The iteration also involves smoothing at two crossings, say $a$ and $b$, which yields \emph{four} terms by (\ref{Tcb}).  In contrast, our construction yields only \emph{one} term corresponding to the single unordered pair $a, b$.  
The two constructions of $\nu_T$ and $\Delta$ (Section~\ref{sec:preliminaries})  differ in their current form, but they share a basic idea and could be connected in a broader framework.  Exploring an approach that relates  Turaev's cobracket to our symmetric setting, as discussed in Section~\ref{sec:Discuss}, would help clarify this relationship.    
\end{remark}
While Turaev's iterated cobracket works in a non-commutative setting, our triple coproduct in $\textup{Sym}^3 (V)$ serves as a commutative counterpart, following   a skein-theoretic principle but adapted to a symmetric tensor context.        
This structural feature underpins the uniqueness result discussed in Section~\ref{sec:Discuss}.    
\section{Discussion}\label{sec:Discuss}
In the previous section, we compared our definition of $\Delta$ with invariants such as the affine index polynomial and noted Turav's cobracket $\nu_T$.  Although these  approaches differ in detail, they share a skein-theoretic basis.  
In this section, we describe their  connection by interpreting our symmetric tensor formulation as part of the classical limit in Turaev's skein quantization.  Here, under the identification $\textup{A}/h\textup{A} \cong \mathbb{Q}[\hbar] \otimes \mathrm{Sym}(V)$, we emphasize  how the $\hbar^2$-term in this expansion corresponds to the image of our triple coproduct. This perspective places our construction as a commutative analogue of Turaev's non-commutative cobracket, extending the framework without replacing it.  

Recalling the definition of $\Delta$ in  (\ref{formula:coproduct}), 
\[
\Delta : \mathcal{D} \to \mathbb{Z}[\mathcal{D}^3]  \otimes W_{LRNP}
\]
maps a knot diagram to a sum of three-component diagrams with over/under information, and requires an explicit proof of invariance under Reidemeister moves as in Section~\ref{sec:ProofM}.  Note that this equivalence, which identifies any triple containing an unknot with zero, ensures well-definedness of the coproduct before applying the intersection form.  

On the other hand, following Professor Tsuji's suggestion, we consider an alternative coproduct:   
\[
\Delta^{\textup{loop}} : \mathcal{D} \to \textup{Sym}^3 (V)  \otimes W_{LRNP}, 
\]
where $V$ is the vector space generated by homotopy classes of free loops on a fixed closed surface.  
This suggestion appears to correspond to a classical limit of the skein algebra of Turaev, specifically the identification $\textup{A}/h\textup{A} \cong \mathbb{Q}[\hbar] \otimes \mathrm{Sym}(V)$, where the $\hbar^2$-coefficient (i.e., the second-order term in the $\hbar$-expansion) corresponds to the image of $\Dloop$ (cf.~\cite[\S 0.3, (0.3.1)--(0.3.2); \S1.3--\S1.4, Theorem~1.4]{Turaev1991}).   
In the framework of Turaev \cite{Turaev1991}, the skein algebra $\textup{A}$ quantizes the Poisson algebra of loops, with $h$ serving as the deformation parameter.  

\noindent{\emph{Convention.}}  In this section, we reserve the symbol $h$ for the deformation parameter controlling non-commutativity in Turaev's skein algebra, and we use the symbol $\hbar$ for Turaev's second parameter, which remains in the coefficient ring after taking the quotient $\textup{A}/h\textup{A}$ \cite{Turaev1991}.  Hence $\textup{A}/h\textup{A}$ is a commutative algebra over $\mathbb{Q}[\hbar]$.   

In \cite{Turaev1991}, the symmetric algebra $S(Z)$ of the Goldman Lie algebra $Z$, generated by  free homotopy classes of loops, appears as the classical limit of $\textup{A}$.  In our setting, the vector space $V$ generated by free homotopy classes plays the role of $Z$, and its symmetric tensor algebra $\textup{Sym}(V)$ corresponds to $S(Z)$ of Turaev.  
Therefore, the identification $\textup{A}/h\textup{A} \cong \mathbb{Q}[\hbar] \otimes \mathrm{Sym}(V)$ in Turaev's  theory aligns with our use of $\mathrm{Sym}(V)$ as the target of the coproduct $\Dloop$.  Furthermore, unlike the standard argument based on Turaev's skein  theory, which establishes well-definedness of the cobracket, our approach demonstrates that, once the target is $\textup{Sym}^3 (V)$, the relations (\ref{eq:RII}) and (\ref{eq:RIII}) reduce the word space  to a one-dimensional quotient.  Consequently, the assignment of weights for smoothing choices becomes uniquely determined.    In this sense, our result extends the known framework based on Turaev's  theory by showing that there is essentially one canonical way to realize the operation in the  symmetric tensor setting.   
This viewpoint was pointed out by Professor Shunshuke  Tsuji in private conversation with NI, one of the authors, and was also discussed in relation to similar consideration by Mr. Hiroki Mizuno in late 2024.          

Since this construction does not require handling over/under information in the target space, the invariance under elementary homotopy moves, similar to `flat' Reidemeister moves, is preserved.  

We specialize this direction to our setting.  For any triple $(C_{ab}^{(1)}, C_{ab}^{(2)}, C_{ab}^{(3)})$ of oriented curves obtained by smoothing a parallel pair of crossings, $a$ and $b$, we define an equivalence relation ``$\sim$'' as follows:   
\[
(C_{ab}^{(1)}, C_{ab}^{(2)}, C_{ab}^{(3)}) \sim 0 \quad \text{if any}~C_{ab}^{(i)}~(i=1, 2, 3)~\text{is an unknot}.  
\]
Here, we write $[C_{ab}^{(1)}, C_{ab}^{(2)}, C_{ab}^{(3)}]$ for the equivalence class under the relation that kills any triple containing an unknot, i.e., we define 
\begin{equation*}
\Delta^{\textup{loop}}(D) := \sum_{a, b : \para}  \varepsilon_a  \varepsilon_b \left[C^{(1)}_{ab}, C^{(2)}_{ab}, C^{(3)}_{ab}\right] \otimes w_{ab}.     
\end{equation*}
\begin{claim}
$\Delta^{\textup{loop}}(D)$ is invariant under stable equivalence, particularly under Reidemeister moves.   
\end{claim}  
\noindent\emph{Sketch of a proof.} First, for the first Reidemeister move, we obtain     
\begin{align*}\label{eq:RoneCB}
\Dloop (D') - \Dloop (D) = \sum_{\substack{a, b : \para \\ a=A~{\textrm{or}~} b=A}}  \varepsilon_a \varepsilon_b  \left[C^{(1)}_{ab}, C^{(2)}_{ab}, C^{(3)}_{ab}\right] \otimes w_{ab} = 0. 
\end{align*} 
Second, we consider the analogue of Lemma~\ref{RIIpara}: 
\begin{equation*}
\left[C^{(1)}_{AB}, C^{(2)}_{AB}, C^{(3)}_{AB}\right] = 0. 
\end{equation*}
Hence, 
\begin{equation*}\label{eqRIIdiffCB}
\begin{split}
&\Dloop (D') - \Dloop (D) = \sum_{\substack{a, b : \para \\ 
|\{ a, b \} \cap \{ A, B \}  |=1}} \varepsilon_a \varepsilon_b  \left[ C^{(1)}_{ab}, C^{(2)}_{ab}, C^{(3)}_{ab} \right] \otimes w_{ab} \\
&=  \sum_{E} \varepsilon_E \varepsilon_A \left[ C^{(1)}_{EA}, C^{(2)}_{EA}, C^{(3)}_{EA} \right] \otimes w_{EA} + \sum_{E} \varepsilon_E \varepsilon_B \left[ C^{(1)}_{EB}, C^{(2)}_{EB}, C^{(3)}_{EB} \right] \otimes w_{EB}\\
&=\sum_{E} \varepsilon_E  \left[ C^{(1)}_{EA}, C^{(2)}_{EA}, C^{(3)}_{EA} \right] \otimes (\varepsilon_A w_{EA} + \varepsilon_B w_{EB}) \quad \left(\because \left[ C^{(1)}_{EA}, C^{(2)}_{EA}, C^{(3)}_{EA} \right]=\left[ C^{(1)}_{EB}, C^{(2)}_{EB}, C^{(3)}_{EB} \right]\right)\\
&=0.  
\end{split} 
\end{equation*}
Third, we consider the third Reidemeister move.   We prove the  invariance by preparing analogous  notations.  Let $A$, $B$, and $C$ be crossings of the triangle involved in  ${\Omega_{3a}}$ and let $D_r$ and $D_{\ell}$ be the diagrams as in  Figure~\ref{IIIaCases}.   For each $\bullet \in \{\ell, r\}$, define  
\begin{align} 
\sum \ast_k (D_{\bullet})^{\Loop} &:=   \sum_{\substack{a, b : \para~{\textrm{in}}~D_{\bullet} \\ 
 |\{ a, b \} \cap \{ A, B, C \}  |= k}} \varepsilon_a \varepsilon_b \left[ C^{(1)}_{ab}, C^{(2)}_{ab}, C^{(3)}_{ab}  \right] \otimes w_{ab}.
 \end{align}    
\begin{itemize}
\item Case of the pair $a, b$ containing exactly \underline{zero} or \underline{one} element of three crossings set $\{ A, B, C \}$.  By definition, $\sum \ast_0 (D_{\ell})^{\Loop}$ $=$ $\sum \ast_0 (D_r)^{\Loop}$.  Using homotopy moves on $\textup{Sym}^3 (V)$ for Figure~\ref{RIIIaDiag}, we also have $\sum \ast_1 (D_{\ell})^{\Loop}$ $=$ $\sum \ast_1 (D_r)^{\Loop}$.   
\item Case of the pair $a, b$ containing exactly \underline{two} elements of three crossings set $\{ A, B, C \}$.  
Before stating this case, let $\Dloop|_{ab} (D)$ $=$ $\varepsilon_a \varepsilon_b  \left[C^{(1)}_{ab}, C^{(2)}_{ab}, C^{(3)}_{ab}\right] \otimes w_{ab}$; then $\displaystyle \Dloop(D) = \sum_{a, b : \para} \Dloop|_{ab} (D)$, as an analogue of Notation~\ref{not:state}.  
Note that any pair of crossings in  $\{ A, B, C \}$ is parallel for one side  $D_{\ell}$ and note also that 
\[\left[C^{(1)}_{AB}, C^{(2)}_{AB}, C^{(3)}_{AB} \right] = \left[C^{(1)}_{AC}, C^{(2)}_{AC}, C^{(3)}_{AC} \right] = \left[C^{(1)}_{BC}, C^{(2)}_{BC}, C^{(3)}_{BC} \right]\] on $\textup{Sym}^3 (V)$.  Then for the invariance for the upper case in Figure~\ref{IIIaCases},  
\begin{align*}
\sum \ast_2 (D_{\ell})^{\Loop} &=\Dloop|_{AB} (D_{\ell})  +  \Dloop|_{BC} (D_{\ell})  + \Dloop|_{AC} (D_{\ell}) \\
&= \left[C^{(1)}_{AB}, C^{(2)}_{AB}, C^{(3)}_{AB} \right] \otimes (RRPP-RRPN-RRNP)\\
&=0.  
\end{align*}
Replacing $RRPP-RRPN-RRNP$ with $LLPP - LLPN - LLNP$, we obtain the invariance for the lower case in Figure~\ref{IIIaCases}.  
Hence,  
\begin{align*}
\Dloop (D_{\ell}) &= \sum \ast_0 (D_{\ell})^{\Loop} +  \sum \ast_1 (D_{\ell})^{\Loop} + \sum \ast_2 (D_{\ell})^{\Loop}  \\
&= \sum \ast_0 (D_{\ell})^{\Loop} +  \sum \ast_1 (D_{\ell})^{\Loop} \\
&=  \sum \ast_0 (D_r)^{\Loop} +  \sum \ast_1 (D_r)^{\Loop} 
= \Dloop (D_r).  
\end{align*}
$\hfill\Box$
\end{itemize}
\begin{corollary}\label{cor:Unique}
Let $D$ be a knot diagram on a surface and let 
\[W = \mathbb{Z} \langle LLPP, LLNP , LLNN, LRPP, LRNP, LRPN, LRNN, RRPP, RRNP, RRNN \rangle,   
\]
where $LLNP$ ($RRNP$,~resp.) is permitted to represent $LLPN$ ($RRPN$,~resp.).   
Then let 
\[
\Delta^{\textup{loop,W}} : \mathcal{D} \to \textup{Sym}^3 (V)  \otimes W; \Delta^{\textup{loop}}(D) := \sum_{a, b : \para}  \varepsilon_a  \varepsilon_b \left[C^{(1)}_{ab}, C^{(2)}_{ab}, C^{(3)}_{ab}\right] \otimes w_{ab}.   
\]
$\Delta^{\textup{loop,W}} (D)$ is invariant under homotopy moves if and only if the elements of $W$ satisfy the relations (\ref{eq:RII}) and (\ref{eq:RIII}).  Moreover, these relations are the only ones  compatible with Reidemeister invariance in the symmetric tensor setting.     
\end{corollary}
\begin{remark}
Although the word space $W_{LRNP}$ becomes effectively one-dimensional after imposing the relations (\ref{eq:RII}) and (\ref{eq:RIII}) derived from Reidemeister moves, we retain it in the construction to explicitly encode the combinatorial data arising from smoothing choices, which determines the structure of the invariant.  This is because this allows us to track how the smoothing choices contribute to the invariant, even if the resulting space is algebraically simple.  
This uniqueness result indicates that, although Turaev's cobracket admits various iterations, the coproduct $\Dloop$ in $\textup{Sym}^3 (V)$ provides a canonical structure under  Reidemeister invariance.  Thus, the construction serves not merely as an analogue but as a natural extension of the existing framework.     
\end{remark}
\section*{Acknowledgements}  
The authors would like to thank Professor Tsuji and Mr.~Mizuno for fruitful discussions regarding the formulation based on loop spaces and its relation to constructions via (quantized) skein algebra of Turaev.  These conversations helped to clarify the theoretical background of the coproduct structure in the current  version.   
The authors also thank Dr.~Keita Nakagane for his careful reading and insightful comments on an earlier version; he identified errors in a previous draft and contributed significantly to improving the exposition.    
The author NI also thanks the referee for valuable comments on earlier versions.  
The work was partially supported by JSPS KAKENHI Grant Numbers  JP20K03604, JPK22K03603, and Toyohashi Tech Project of Collaboration with KOSEN Grant Number 2309.  
A preliminary version of this work was initiated while the authors were affiliated with National Institute of Technology, Ibaraki College, Japan, and partially based on an undergraduate project by TK, one of the authors, conducted under the supervision of Professor Mariko Okude.    
While the present version does not incorporate ideas such as an infinite sequence, that project served as motivation for this study, and Professor Okude's encouragement is gratefully acknowledged.
\bibliographystyle{plain}
\bibliography{Ref}
\end{document}